\documentclass[]{ajour}
\usepackage{amsmath, amssymb, amsthm, conformal, mystyle}
\usepackage{hyperref}
\hypersetup{colorlinks=false, pageanchor=true, breaklinks=true}

\renewcommand\^[1]{\widehat{#1}}
\renewcommand{\theenumi}{{\rm (\alph{enumi})}}

\def\c{{\mathsf c}}

\begin{document}

%


\authorrunninghead{Michael Roitman}
\titlerunninghead{Combinatorics of free vertex algebras}



\title{Combinatorics of free vertex algebras}

\author{Michael Roitman}

\affil{M.S.R.I.\\
1000 Centennial Dr.\\
Berkeley, CA 94720}
\email{roitman@msri.org}




\begin{article}

\vfill
\section*{Introduction}
\hypertarget{intro}{}
This paper illustrates the combinatorial approach to vertex algebra
--- study  of vertex algebras  presented by generators and
relations. A necessary ingredient of this method is the notion of
free vertex algebra. Borcherds \cite{bor} was the first to note 
that free vertex algebras do not exist in general. The reason for this
is that vertex algebras do not
form a variety of algebras, as defined in
e.g. \cite{cohn}, because the \hyperlink{V4}{locality axiom} (see \sec{def}
below)  is not an identity. However, a certain subcategory of vertex
algebras, obtained by restricting the \hyperlink{order}{order of locality} of generators,
has a universal object, which we call the free
vertex algebra corresponding to the given locality bound. 
In \cite{freecv} some free vertex algebras were
constructed and in certain special cases their linear bases were
found. In this paper 
we generalize the construction of \cite{freecv} and find linear bases
of an arbitrary free vertex algebra. 

It turns out that free vertex algebras are closely related to the
vertex algebras corresponding to integer lattices. The latter algebras play
a very important role in different areas of mathematics and
physics. They were extensively studied in e.g. \cite{dong,dl,flm,kac2,lixu}.
Here we explore the relation between free vertex algebras and lattice
vertex algebras in much detail. These results comply with the use of
the word ``free'' in physical literature refering to some elements of 
lattice vertex algebras, like in 
``free field'', ``free bozon'' or ``free fermion''.  

Among other things, we find a nice
presentation of lattice vertex algebras in terms of generators and
relations, thus giving an alternative construction of these algebras
without using vertex operators. We remark that our construction works
in a very general setting; we do not assume the lattice to be 
positive definite, neither  non-degenerate, nor of a finite rank. 

\vfill

\subsubsection*{Organization of the manuscript}
We start with reviewing some basic definitions of
\hyperlink{V1}{vertex} and \hyperlink{conformal}{\hbox{conformal}} 
algebra in \sec{fields}--\sec{conformal}. The reader may consult Kac's
book \cite{kac2} for more 
details. Then in \sec{lattice} we recall the construction of the
vertex algebra $V_\Lambda$ corresponding to an integer lattice
$\Lambda$, the details can be again found in \cite{kac2}. 

In \sec{free} we construct the free vertex algebra $\goth F_N(\cal B)$
generated by a set $\cal B$ so that the \hyperlink{order}{order of
locality} $N(a,b)$ of a pair
of generators $a,b\in\cal B$ is given by an arbitrary symmetric 
integer-valued function $N:\cal B\times\cal B\to \Z$. \thm{basis} describes a set 
$\cal T\subset\goth F_N(\cal B)$ and claims that $\cal T$ is a linear 
basis of $\goth F_N(\cal B)$. As a
corollary, we find the
dimensions of homogeneous components of $\goth F_N(\cal B)$. 
In \sec{span} we will prove that $\cal T$ spans
$\goth F_N(\cal B)$, and  in   \sec{linind} we finish the proof of 
\thm{basis} by showing that $\cal T$ is linearly independent.

In \sec{embedding} we construct a vertex algebra homomorphism
$\f:\goth F_N(\cal B)\to V_\Lambda$ from the free vertex algebra 
$\goth F_N(\cal B)$ to the vertex algebra corresponding to the
lattice $\Lambda = \Z[\cal B]$. The integer form on $\Lambda$ is
defined by $(a|b) = -N(a,b)$ for $a,b\in\cal B$. \thm{embedding} then states
that $\f$ is injective. We will prove this theorem in \sec{linind}.

In \sec{dong}, using \thm{embedding}, 
we prove a quantitative version of \hyperlink{dong}{Dong's lemma}. 
In \sec{locfun} we apply this
lemma to settle a question raised in \cite{universal}: we prove that the
\hyperlink{locfun}{locality function} of a free
\hyperlink{conformal}{conformal} \hyperlink{conformal}{algebra} has 
quadratic growth.  

In \sec{confact} we study homogeneous \hyperlink{confder}{conformal
derivations}  of free
vertex algebras. It turns out that a particularly interesting case is
when the algebra is generated by a single element $a$ such that
$N(a,a)=-1$. In this case $\goth F_N\big(\{a\}\big)$ is embedded
into the fermionic vertex algebra $V_\Z$. We prove that a
homogeneous component of 
$\f\big(\goth F_N(\{a\})\big)\subset V_\Z$ is an 
irreducible lowest weight module over certain conformal algebra 
$\^{\goth W}\subset V_0$, such that the\break
\hyperlink{coeff}{coefficient algebra} of $\^{\goth W}$ is a central
extension of the Lie 
algebra of differential operators on the circle, see \cite{fkrw,kac2,bozfer}.

Finally, in \sec{present} we find a presentation of lattice vertex
algebras in terms of generators and relations, see \thm{present}. 
It turns out that the
required relations are rather minimal. Our proof is completely
combinatorial. The first step is to determine the structure of 
the vertex algebra in question as a module over the
\hyperlink{heisen}{Heisenberg algebra}. 

\vfill

\subsubsection*{Notations}
All algebras and linear spaces are over a field $\Bbbk$ of
characteristic 0. Here are some shortcuts used throughout the paper: 
$\set{n\in\Z}{n\ge 0}=\Z_+$, \ \ 
$\frac 1{k!} n(n-1)\cdots (n-k+1)=\binom nk$, \ \ 
$\frac 1{k!}\,D^k=D^{(k)}$, \ \ $Z(L)$ is the center of a Lie algebra $L$.

\vfill

\hypertarget{ack}{}
\subsubsection*{Acknowledgements}
The core of this work was done while I was visiting the
Fields Institute. 
I am very grateful to organizers of the Program on Lie Theory
for inviting me there. I especially thank Stephan
Berman, Yuly Billig and Chongying Dong for helpful discussions.  

\section{Fields and locality}\label{sec:fields}
Let $V=V\even\oplus V\odd$ be a vector superspace over $\Bbbk$. 
Take a formal variable $z$ and 
consider the space 
$\F(V)=\F(V)\even\oplus\F(V)\odd \subset \op{End}\big(V,V((z))\big)$
of {\it fields} on $V$, given by
$$
\F(V)^p = \left\{ \left.
\sum_{n\in \mathbb Z}a(n)\,z^{-n-1}\ \right|
p\big(a(n)\big)=p, \ 
\forall v \in V, \  a(n) v = 0 \text{  for } n \gg 0
\right\},
$$
Here $p(x)\in \Z/2\Z$ is the parity of $x$.
Denote by  $\1 \in \F(V)\even$ the identity operator, 
such that $ \1 (-1) = \op{Id}_V$,  all other coefficients are 0. 
 
Let $\imath_{w,z}(w-z)^n$ and  $\imath_{z,w}(w-z)^n$ be two
different expansions of $(w-z)^n$ into formal power series in the
variables $w$ and $z$:
\begin{equation*}
\begin{split}
\imath_{w,z}(w-z)^n &=\sum_{i\ge0}(-1)^{n+i} \binom ni\, w^{n-i}z^i 
\in \Bbbk[[w,w\inv,z]], \\
\imath_{z,w}(w-z)^n &=\sum_{i\ge0}(-1)^i \binom ni\, w^i z^{n-i}
\in \Bbbk[[w,z,z\inv]].
\end{split}
\end{equation*}
Of course, if $n\ge 0$ then $\imath_{w,z}(w-z)^n=\imath_{z,w}(w-z)^n$.

Let $a, b \in \F(V)$. We say that  $a$ is {\it
local} to $b$ if there is some $N\in\Z$ such
that 
\begin{equation}\label{fl:locality}
a(w)b(z)\,\imath_{w,z}(w-z)^N - (-1)^{p(a)p(b)}
b(z)a(w)\,\imath_{z,w}(w-z)^N = 0.
\end{equation}
The minimal $N=N(a,b)$ with this property is called the 
\hypertarget{order}{{\it order of locality}} of $a$ and $b$. 
In terms of the coefficients locality means that the
following identities hold for all $m,n\in\Z$:
\begin{equation}\label{fl:locoef}
\begin{split}
\sum_{s\ge0}&(-1)^s \binom Ns a(m-s)b(n+s)\\
&-(-1)^{p(a)p(b)}
\sum_{s\le N}(-1)^s \binom N{N-s} b(n+s)a(m-s)=0.
\end{split}
\end{equation}

\noindent
Usually the locality is defined only for a nonnegative order $N$
 \cite{dl,kac2}, but we
will need this property in a bigger generality.

For an integer $n$ define a bilinear product 
$\ensquare n :\F(V)\otimes \F(V) \to \F(V)$ by 
\begin{equation}\label{fl:prod}
\begin{split}
\big(a \ensquare{n} b\big)(z) = 
\op{Res}_w \biggl(
a(w)&b(z)\,\imath_{w,z}(w-z)^n\\ 
&-(-1)^{p(a)p(b)}
b(z)a(w)\,\imath_{z,w}(w-z)^n
\biggr).
\end{split}
\end{equation}
Here $\op{Res}_w$ stands for the coefficient of $w\inv$. Explicitly
this means
\begin{equation}\label{fl:prodcoef}
\begin{split}
\big(a \ensquare{n} b\big)(m) &= 
\sum_{s\ge0}(-1)^s \binom ns a(n-s)b(m+s)\\
&-(-1)^{p(a)p(b)}\sum_{s\le n}(-1)^s \binom n{n-s} b(m+s)a(n-s).
\end{split}
\end{equation}
Clearly, if $a$ and $b$ are local of order $N$ then 
$a \ensquare n b = 0$ for $n\ge N$. 

If $n=-1$, then $a \ensquare{-1}b = \:ab\:=
a_-b + ba_+$ is the so-called {\it normally ordered
product}, where 
$$
a_{\pm}(z) = \smash{\sum_{n\genfrac{}{}{0pt}{}{\ge}{<} 0}} a(n)\,z^{-n-1}.  
$$
In general, if $n<0$ then 
$$
\alpha \ensquare n \beta =  \:\big(D^{(-n-1)}\alpha\big)\beta\:, 
$$
where $D\alpha(z) = \frac d{dz}\alpha(z) = \alpha(z) \ensquare{-2}\1$.

We remark that the following identities hold:
$$
(D a)\ensquare n b = -n\, a \ensquare{n-1} b,\quad
a\ensquare n (Db) = n\, a \ensquare{n-1} b + D\big(a\ensquare n b\big), 
$$
in particular, $D$ is a derivation of the products $\ensquare n$.

The Dong's lemma \cite{dl,kac2,li} states that if  
$a,\,b,\,c \in F(V)$ are three pairwise local fields,
then for every $n \in \Z$ the fields $a \ensquare n b$ and $c$ are 
local as well. In \sec{dong} we will prove a quantitative version of this
lemma.   

A subspace $\goth A\subset \F(V)$ such that 
\begin{itemize}
\item[(i)]
any two fields $a,b\in\goth A$ are local,
\item[(ii)]
$\goth A\ensquare n \goth A \subseteq \goth A$ for all $n\in\Z$ and 
\item[(iii)]
$\1\in\goth A$
\end{itemize}
is called a vertex
algebra. By the \hypertarget{dong}{Dong's lemma}, 
given a collection $\cal S\subset \F(V)$
of pairwise local fields, the closure of $\cal S\cup\{\1\}$
under the products is a vertex algebra in $\F(V)$. 

\section{Axiomatic definition of vertex algebras}\label{sec:def}
Vertex algebras can also be defined axiomatically as follows
\cite{kac2}. Let $\goth A=\goth A\even\oplus\goth A\odd$ be a linear 
superspace endowed with a sequence of even bilinear
operations $\ensquare n:\goth A\otimes \goth A\to \goth A,\ n\in \Z$,
and a distinguished element $\1\in \goth A\even$. Let  $D:\goth A\to
\goth A$ be an even linear map given by $Da = a\ensquare{-2}\1$. 
Consider the {\it left regular action map} $Y:\goth A\to \goth A[[z,z\inv]]$ 
defined by $Y(a)(z) = \sum_{n\in\Z}(a\ensquare n \,\cdot\,)\,z^{-n-1}$.
Then $\goth A$ is a vertex algebra if it satisfies  
the following conditions for any $a,b\in \goth A$ and $n\in \Z$:

\begin{itemize}
\item[V1.]\hypertarget{V1}{}
$a\ensquare n b = 0$ for $n\gg0$.
\item[V2.]\hypertarget{V2}{}
$\1\ensquare n a = \delta_{n,-1}a,\qquad
a\ensquare n \1 = \begin{cases}
0&\text{if}\quad n\ge0,\\
D^{(-n-1)} a&\text{if}\quad n<0.
\end{cases}$
\item[V3.]\hypertarget{V3}{}
$(D a)\ensquare n b = -n\, a \ensquare{n-1} b,\quad
a\ensquare n (Db) = n\, a \ensquare{n-1} b + D\big(a\ensquare n b\big)$.
\item[V4.]\hypertarget{V4}{}
The series $Y(a), Y(b)\in \goth A[[z,z\inv]]$ are local.
\end{itemize}
This is not a minimal set of axioms, in fact conditions V2 and V3 can
be weakened. Note that V3 implies that $Y(Da) = \frac d{dz}\, Y(a)$.
For other axiomatic definitions of vertex algebras see \cite{bor,dl,fhl,flm}.

A homomorphism  of two vertex algebras 
$\goth A$ and $\goth B$ is a map $\phi:\goth A\to\goth B$ such that 
$\phi(\1)=\1$ and $\phi(a\ensquare n b)=\phi(a)\ensquare n\phi(b)$ for
all $n\in \Z$, in particular $\phi(Da) = \frac d{dz}\, \phi(a)$. 
A module over a vertex algebra $\goth A$ is a linear
space $V$ such that there is a vertex algebra homomorphism 
$\phi:\goth A\to \F(V)$ of $\goth A$ into a vertex subalgebra of
$\F(V)$. Similarly one defines the notions of ideals, quotients,
isomorphisms, etc.

One of the most important properties of vertex algebras is that the
left regular action map is in fact  an isomorphism of a vertex algebra
$\goth A$ with a vertex subalgebra  $Y(\goth A)\subset\F(\goth A)$. 
In other words, the left regular action map is a faithful
representation of $\goth A$ on itself. In view of \fl{prodcoef}
this implies the following identity:
\begin{equation}\label{fl:assoc}
\begin{split}
(a \ensquare{n} b)\ensquare m c &= 
\sum_{s\ge0}(-1)^s \binom ns a\ensquare{n-s}(b\ensquare{m+s}c)\\
&-(-1)^{p(a)p(b)}\sum_{s\le n}(-1)^s \binom n{n-s}
b\ensquare{m+s}(a\ensquare{n-s}c)
\end{split}
\end{equation}
for all $m,n\in\Z$, \ $a,b,c\in\goth A$. An equivalent form of this
identity is 
\begin{equation}\label{fl:comm}
a\ensquare m \big(b\ensquare n c\big)
- b\ensquare n \big(a\ensquare m c\big) = 
\sum_{s\ge 0}\binom ms \big(a\ensquare s b\big)\ensquare{m+n-s}c.
\end{equation}
Another important idenity is the following so-called quasisymmetry
identity:
\begin{equation}\label{fl:qs}
a\ensquare n b = - \sum_{i\ge0} (-1)^{n+i} D^{(i)}
(b\ensquare{n+i}a).
\end{equation}

There is a couple of additional axioms which are often imposed on vertex
algebras, see \cite{dl,fhl,flm,li}. A vertex algebra $\goth A$ is
called a {\it vertex operator algebra} or a {\it conformal vertex
algebra} if
\begin{itemize}
\item[V5.]\hypertarget{V5}{}
$\goth A =\bigoplus_{i\in\frac12\Z} \goth A_i$ is graded so that 
$\goth A\even = \sum_{i\in \Z}\goth A_i$, \ 
$\goth A\odd = \sum_{i\in \Z+\frac12}\goth A_i$, \ 
$\1\in \goth A_0$ and 
$\goth A_i \ensquare n \goth A_j \subset \goth A_{i+j-n-1}$.
\item[V6.]\hypertarget{V6}{}
There exists an element $\omega \in \goth A_2$  such that
$N(\omega, \omega)=4$,\ 
$\omega \ensquare 0 a = Da$, \ 
$\omega \ensquare 1 a = (\deg a)\, a$ for all homogeneous $a\in\goth A$,\ 
$\omega \ensquare 3 \omega = 0$ and 
$\omega \ensquare 4 \omega = c\1$ for some $c\in \Bbbk$.
\end{itemize}
The axiom V6 implies that $\omega$ generates the the \hyperlink{Vir}{Virasoro
conformal algebra} $\goth{Vir} \subset \goth A$ 
(see \sec{conformal} below).
The number $c$ is called the {\it conformal charge} of $\goth A$.

We will call an expression of the form 
$a_1 \ensquare{n_1} \ldots \ensquare{n_{l-1}}a_l$, \ 
$a_i\in\goth A$, \ $n_i\in\Z$, (with an arbitrary
order of parentheses) a {\it vertex monomial} of length $l$.  
Using \fl{assoc} every vertex monomial can be expressed
as a linear combination of {\it right normed} monomials of the form
$a_1 \ensquare{n_1}( \,\cdots a_{l-2}  \ensquare{n_{l-2}}
(a_{l-1} \ensquare{n_{l-1}}a_l)\mbox{$\cdot\!\cdot\!\cdot$})$.
We will call a monomial $w$ \hypertarget{minimal}{{\it minimal}} 
if every its submonomial $u\ensquare n v$ has the
maximal possible value of $n$, i.e. $n=N(u,v)-1$.

\section{Conformal algebras}\label{sec:conformal}
A \hypertarget{conformal}{conformal algebra} is a more general algebraic structure than 
a vertex algebra, see \cite{kac2}. The former is obtained from the latter by
forgetting about $\1$ and the products $\ensquare n$ for $n<0$. 
The identities \hyperlink{V1}{V1}, \hyperlink{V3}{V3}, \fl{assoc} and
\fl{qs} hold in a conformal 
algebra for $m,n\ge 0$.

Note that if $N\ge 0$ and $n\ge 0$, then \fl{locality} and \fl{prod}
simplify respectively  to 
\begin{gather}\label{fl:posloc}
\ad{a(w)}{b(z)}\,(w-z)^N = 0,\\
\label{fl:posprod}
\big(a \ensquare{n} b\big)(z) = 
\op{Res}_w \ad{a(w)}{b(z)}\,(w-z)^n.
\end{gather}
\noindent
This means that the locality for $N\ge 0$ and the
products $\ensquare n$ for $n\ge 0$ both make sense for series $a,b \in L[[z^{\pm1}]]$
with coefficients in a Lie superalgebra $L$. Similar to the case of
vertex algebras, a space $\goth L \subset  L[[z^{\pm1}]]$ of pairwise
local series,
such that $D\goth L \subseteq \goth L$ and $\goth L \ensquare n \goth L
\subseteq \goth L$ for $n\ge 0$, is a conformal algebra and all
conformal algebras can be obtained in this way. 
The \hyperlink{dong}{Dong's}  \hyperlink{dong}{lemma} holds
for these series as well.

Moreover, for a given conformal algebra $\goth L$ one can construct a
Lie superalgebra $L=\cff \goth L$, which is called the 
\hypertarget{coeff}{{\it coefficient algebra}}
of $\goth L$, together with a conformal algebra homomorphism $\goth L \to  L[[z^{\pm1}]]$
which is universal among all representations of $\goth L$ by formal
series with coefficients in a Lie superalgebra. As a linear space
$L$ is equal to $\goth L \otimes \Bbbk[t^{\pm1}]$ modulo the linear
subspace generated by all the expressions $(Da)(n) + n a(n-1)$ for 
$a\in \goth L$ and $n\in \Z$. We denote here $a(n) = a\otimes t^n$.

\subsubsection*{Examples}
Let $\goth g$ be a Lie algebra with a symmetric invariant bilinear
form $(\,\cdot\,|\,\cdot\,)$. Consider the corresponding affine Lie
algebra $\tl{\goth g} = \goth g\otimes \Bbbk[t^{\pm1}] \oplus \Bbbk\c$
with the brackets given by 
$$
\ad{a(m)}{b(n)} = \ad ab(m+n) + \delta_{m,-n}\, m (a|b)\,\c, \qquad 
\c\in Z(\tl{\goth g}).
$$
This Lie algebra is the coefficient algebra
of the conformal algebra $\goth G \subset \tl{\goth g}[[z^{\pm1}]]$,
generated by  the series $\tl a= \sum_n a(n)\,z^{-n-1}$ for $a\in \goth g$ and
$\c=\c(-1)$ so that $D\c=0$,\ $N(\tl a, \tl b)=1$ and 
$$
\tl a\ensquare{0} \tl b = \tl {\ad ab},\qquad
\tl a\ensquare{1} \tl b = (a|b)\,\c.
$$
In the case when $\goth g$ is an abelian Lie algebra, the
corresponding affine algebra $\tl{\goth g}$ is a
\hypertarget{heisen}{Heisenberg algebra}, and $\goth G$ is a Heisenberg 
conformal algebra.

Another example of a conformal algebra is obtained from the Virasoro
Lie algebra $Vir$, which is spanned by the elements $L_n$, \ $n\in\Z$ and $\c$
with the brackets given by 
$$
\ad{L_m}{L_n} = (m-n) L_{m+n} +
\delta_{m,-n}\binom {m+1}3 \c, \qquad \c\in Z(Vir). 
$$
It is the 
coefficient algebra of the conformal algebra 
\hypertarget{Vir}{$\goth{Vir}\subset Vir[[z^{\pm1}]]$} generated by the series 
$\omega = \sum_{n\in\Z}L_{n-1}\,z^{-n-1}$ and $\c$. We have $D\c=0$, \
$N(\omega, \omega)=4$ and 
$$
\omega\ensquare0 \omega = D\omega, \quad
\omega\ensquare1 \omega = 2\omega, \quad
\omega\ensquare2 \omega = 0, \quad
\omega\ensquare3 \omega = \c.
$$

Another example of a conformal algebra will be given in \sec{confact}.

\subsubsection*{Conformal operators}
Consider first a $\Bbbk[D]$-module $V$. A {\it conformal operator} 
\cite{ck} on
$V$ is a series $\alpha = \sum_{n\in\Z_+}\alpha(n)\,z^{-n-1} \in 
z\inv(gl\, V)[[z\inv]]$ such that $\alpha(n)v = 0$ if $n\gg0$ 
for every fixed $v\in V$ and $\ad D{\alpha} = \frac d{dz} \alpha$.  
For $v\in V$ call $N(\alpha,v) = \min\set{n\in\Z_+}{\alpha(m)v=0\
\forall\,m\ge n}$ the order of locality of $\alpha$ and $v$.
Denote by $cgl\, V$ the space of all conformal
operators on $V$.

A pair $\alpha, \beta \in cgl\,V$ is said to be local of
order $N\in\Z_+$ if \fl{posloc} holds. Also the products $\alpha \ensquare n
\beta\in  cgl\,\goth A$ are defined for $n\ge 0$ by the formula
\fl{posprod}.  The \hyperlink{dong}{Dong's lemma}
 also holds for conformal 
operators. A subspace $\goth L \subset cgl\,V$ of pairwise local
conformal operators, closed under all the products $\ensquare n$ for
$n\in\Z_+$, is a conformal algebra. A module over a conformal algebra
$\goth L$ is a $\Bbbk[D]$-module $V$ with a conformal algebra homomorphism 
$\goth L \to cgl\,V$.

It is  well-known \cite{ck} that if $v\in V$ is such that
$p(D)v=0$ for some polynomial $p\in \Bbbk[D]$ then $\alpha(n)v=0$ for
every $\alpha\in cgl\, V$ and $n\in \Z_+$. In particular, if $\goth A$
is a vertex algebra, then $\alpha(n)\1=0$ for every $\alpha \in
cgl\, \goth A$.

Let $\goth A$ be a conformal (in particular, vertex) algebra. 
A \hypertarget{confder}{{\it conformal derivation}} \cite{ck,universal} is a conformal
operator $\alpha \in cgl\, \goth A$ such that 
$$
\alpha(m) (a\ensquare n b) = a\ensquare n \big(\alpha(m)b\big) + 
\sum_{s\ge 0}\binom ms \big(\alpha(s)a\big)\ensquare{m+n-s}b
$$
for all $m\in\Z_+,\, n\in \Z$. For example, if $a \in \goth A$, then 
$$
Y(a)_+=\sum_{n\in\Z_+} Y(a)(n)\,z^{-n-1}
$$ 
is a conformal derivation of
$\goth A$. Denote by $cder\, \goth A\subset cgl\,\goth A$ 
the space of all conformal derivations of $\goth A$.
It is not difficult to show that if 
$\alpha, \beta \in cder\, \goth A$, then   $\alpha \ensquare n
\beta \in cder\, \goth A$.

Assume that $\goth A$ is generated by a set $\cal B$. Then a conformal
derivation is uniquely defined by its action on $\cal B$. It is easy to show
that if for $\alpha, \beta \in cder\, \goth A$ 
the orders of locality $N(\alpha, \cal B)$, \   
$N(\beta, \cal B)$ are uniformly bounded on $\cal B$, then $\alpha$ and
$\beta$ are local.

\section{Lattice vertex algebras}\label{sec:lattice}
In this section we give a very important example of vertex
algebras --- the algebra $V_\Lambda$ corresponding to an integer
lattice $\Lambda$.  We mostly follow 
\cite{kac2}, see also \cite{dong,dl,flm}. Note, however, that we do
not assume that $\Lambda$ is of a finite rank.

Denote the  integer-valued symmetric bilinear form on $\Lambda$ by 
$(\,\cdot\,|\,\cdot\,)$. Let us extend the form to 
$\goth h = \Lambda \otimes_\Z \Bbbk$. 
Let $H=\tl {\goth h} = \goth h\otimes \Bbbk[t,t\inv]\oplus \Bbbk\c$
be the corresponding \hyperlink{heisen}{Heisenberg algebra}, see \sec{conformal}. 

For every $a \in \Lambda$ consider the canonical relation
level 1 representation $V_a$  of $H$, that is, 
the irreducible Verma $H$-module generated by the vacuum
vector $v_a$ such that 
$h(n)v_a = 0$ for $n>0$, \  $h(0)=(h|a)\,\op{Id}$, and 
$\c=\op{Id}$, see \cite{kac1}. It is well-known (see e.g. \cite{kac2})
that the module $V_0$ has a structure of vertex algebra, such that 
$v_0 = \1$. The map $\tl h \mapsto h(-1)v_0$ is an injective conformal
algebra homomorphism of the
conformal Heisenberg algebra $\goth H$
to $V_0$. From now on we will identify $\goth H$ with its image in $V_0$.    
In fact one can show that $V_0$ is a unique vertex algebra generated
by $\goth H$. 

For $a\neq 0$, 
the map $\tl h\mapsto \sum_{n\in\Z} h(n)\,z^{-n-1}\in\F(V_a)$ is
a conformal algebra homomorphism of $\goth H$ into $\F(V_a)$,
which extends to the vertex algebra
homomorphism $V_0 \to \F(V_a)$, so that $V_a$ becomes a
module over the vertex algebra $V_0$.
Define the Fock space 
$$
V_\Lambda=\bigoplus_{a\in\Lambda} V_a \cong
V_0\otimes\Bbbk[\Lambda].
$$ 

Let $\e:\Lambda\times \Lambda \to \{\pm 1\}$ be a
bimultiplicative map such that 
\begin{equation}\label{fl:eps}
\e(a,b) = (-1)^{(a|a)(b|b)} 
(-1)^{(a|b)}\e(b,a).
\end{equation}
for any $a,b\in\Lambda$. We remark that it is enough to check 
the identity \fl{eps} only when $a$ and $b$ belong to some 
$\Z$-basis of $\Lambda$; then \fl{eps} will follow for general
$a, b$ by bimultiplicativity. The cocycle $\e$ defines an element of 
$H^2\big(\Lambda, \{\pm1\}\big)$. 

The main result \cite{dong,dl,flm,kac2} is that the vertex algebra 
structure on $V_0$ can be uniquely extended to  $V_\Lambda$ such that
$\tl h \ensquare n v = h(n)v$ for any  
$h\in \goth h$ and $v\in V_\Lambda$. The locality of a pair
$v_a,v_b\in V$ of vacuum vectors is
$N(v_a,v_b)=-(a|b)$ and the products are
\begin{equation}\label{fl:vanvb}
v_a \ensquare{-(a|b)-k-1} v_b = 
\e(a,b)\, \big(D-b(-1)\big)^{(k)} v_{a+b},
\qquad k\ge 0.
\end{equation}
In particular, 
$$
v_a \ensquare{-(a|b)-1} v_b = \e(a,b)\, v_{a+b}, \qquad 
v_a \ensquare{-(a|b)-2} v_b = \e(a,b)\, a(-1)v_{a+b}.
$$ 
Taking $b=-a$ we get
$$
v_a \ensquare{(a|a)-1}v_{-a} = \e(a,a)\,\1, \qquad 
v_a \ensquare{(a|a)-2}v_{-a} = \e(a,a)\,\tl a.
$$
It follows that the vacuum
vectors $v_{\pm a}$, for $a$ running over an integer basis of $\Lambda$,
generate the vertex algebra $V_\Lambda$, and \fl{vanvb}  defines the
vertex algebra structure on $V_\Lambda$ uniquely.
The vertex algebra $V$ is simple if and only if the form
$(\,\cdot\,|\,\cdot\,)$ is non-degenerate.


The standard way to construct vertex algebra $V_\Lambda$ is by
showing that the vacuum vectors $v_a$ act under the left regular
map $Y:V_\Lambda \to \F(V_\Lambda)$ by the so-called vertex operators.
In \sec{present} we will construct $V_\Lambda$ by a different
method.

Besides the grading by the lattice $\Lambda$, the vertex algebra 
$V_\Lambda$ has a grading by $\frac 12\Z$, so that 
$\deg v_a = \frac12 (a|a),\ \deg \tl h = 1$ for every
$h\in \goth h,\ \deg \ensquare n = -n-1$ and $\deg D = 1$, in
particular the
axiom \hyperlink{V5}{V5} always holds. We have
decomposition 
$$
V_a = 
\bigoplus_{d\in \frac12 (a|a) + \Z_+} V_{a,d},\qquad 
V_{a,d}\ensquare n V_{a',d'} \subset V_{a+a', d+d'-n-1}.
$$
We will refer to the grading by $\Lambda$ as the grading by weights,
and write $\op{wt}u = a\in \Lambda$ for $u\in V_a$.

Consider the case when $\Lambda$ is non-degenerate and of a finite rank
$l$.  Let $a_1,\ldots,a_l$ and
$b_1,\ldots,b_l$ 
be dual bases of $\goth h$, i.e. such that
$(a_i|b_j)=\delta_{ij}$. Then the element 
$\omega = \frac 12\sum_{i=1}^l a_i \ensquare{-1} b_i \in V_0$
generates a copy of the Virasoro conformal algebra $\goth{Vir}$
so that  \hyperlink{V6}{V6} is satisfied, and hence   
$V_\Lambda$ becomes a vertex operator algebra.

\begin{Rem}
The correspondence $\Lambda \mapsto V_\Lambda$ is not a functor,
because there is a certain degree of freedom in choosing the cocycle
$\e$, satisfying \fl{eps}. However, all the resulting vertex algebras
$V_\Lambda$ are isomorphic.
\end{Rem}
\section{Free vertex algebras}\label{sec:free}
In this section we construct another example of vertex algebras --- a
free vertex algebra. Take some set of generators $\cal B$. We will
assume that $\cal B$ is linearly ordered. Suppose we have a symmetric
function $N: \mathcal B \times  \mathcal B \to \Z$.
Consider the category $\goth{Ver}_N(\cal B)$ of  
vertex algebras 
generated by the set $\mathcal B$ such that in any vertex algebra
$\goth A \in \goth{Ver}_N(\cal B)$ one has
$a\ensquare n b = 0$ for any $a,b \in \cal B$ whenever $n \ge N(a,b)$.
We set the parity of $a\in \cal B$ to be  $p(a) \equiv N(a,a)\mod2$.
The morphisms of $\goth{Ver}_N(\cal B)$ are, naturally, vertex
algebra homomorphisms $\phi:\goth A_1 \to \goth A_2$ such that 
$\phi(a) = a$ for any 
$a \in \mathcal B$. Then it is easy to see that $\goth{Ver}_N(\cal B)$
has a unique universal object $\goth F_N(\cal B)$, 
called {\it the free vertex algebra
generated by $\cal B$ with respect to locality bound $N$}. We
construct $\goth F_N(\cal B)$ explicitly below.

\begin{Rem}
It follows easily from the quazisymmetry identity \fl{qs} that if
$a\in \goth A\even$ is  
an even element of a vertex algebra $\goth A$, then $N(a,a)\in 2\Z$.
\end{Rem}

Let $\goth A$ be a vertex algebra generated by a set $\cal B$. For
$a,b\in\cal B$ let
$$
N(a,b) = \min\bigset{n\in\Z}{a\ensquare m b =0 \ \ \forall\, m\ge n}
$$ 
be the order of locality. Then there is a surjective homomorphism 
$\psi:\goth F_N(\cal B)\to \goth A$. Let 
$\cal R \subset \goth F_N(\cal B)$ 
be a set of generators of $\op{Ker}\psi$. Then we say that $\goth A$
is presented by generators $\cal B$, locality bound 
$N:\cal B\times\cal B\to \Z$ and relations
$\cal R$.

We construct the free vertex algebra $\goth F_N(\cal B)$ as follows. Consider the set  
$\cal X=\set{a(n)}{a\in\cal B,\ n\in\Z}$ and let $\Bbbk[\cal X]$ be the free
associative algebra generated by $\cal X$. A module $M$ over $\Bbbk[\cal X]$ is
called {\it restricted} if for any $a\in \cal B$ and $x\in M$ we have 
$a(n)x=0$ for $n\gg0$. Let $\ol{\Bbbk[\cal X]}$ be the completion of $\Bbbk[\cal X]$
with respect to the topology in which any series $\sum_i p_i$, \ 
$p_i\in \Bbbk[\cal X]$, that makes sense as an operator on every restricted
module, converge. Let $I\subset \ol{\Bbbk[\cal X]}$ be the ideal generated
by the locality relations \fl{locoef} for all $a,b\in\cal B$, \ $N=N(a,b)$.
Denote $U=\ol{\Bbbk[\cal X]}/I$. Let 
$\cal X_+ = \set{a(n)}{a\in\cal B,\ n\ge 0}\subset \cal X$ and consider the
quotient  $U/U\cal X_+$. This is a restricted left $U$-module, therefore
any $a\in \cal B$ corresponds to a field 
$\sum_{n\in\Z} a(n)\,z^{-n-1}\in \F(U/U\cal X_+)$. These fields are local
by the construction and hence they generate a vertex algebra 
$\goth F_N(\cal B)\subset \F(U/U\cal X_+)$, which clearly satisfies all the
criteria for being the desired free vertex algebra. By the Goddard-Kac 
existence theorem \cite[Theorem 4.5]{kac2}, the map 
$\phi \mapsto \phi(-1)1$ is a one-to-one correspondence between 
$\goth F_N(\cal B)$ and $U/U\cal X_+$, so from now on we will identify the two.

In the special case when  $N(a,b)\ge0$ the locality identities \fl{locoef} are
finite sums of commutators, hence instead of $U$ one can consider an
enveloping algebra of a certain Lie algebra. In this case the
construction of $\goth F_N(\cal B)$  was done in
\cite{freecv}. It was also proved there that if $N\equiv N(a,b)\in2\Z_+$ is
a non-negative even constant, then a linear basis of  $\goth F_N(\cal B)$ 
is given by all right normed vertex monomials 
\begin{equation}\label{fl:monom}
a_1 \ensquare{m_1} (a_2  \ensquare{m_2}
\cdots (a_k \ensquare{m_k}\1)\mbox{$\cdot\!\cdot\!\cdot$}),\qquad 
a_i\in\cal B, \ m_i\in\Z, \ m_k <0,
\end{equation}
such that
\begin{equation*}
m_i - m_{i+1} \le 
\begin{cases}
N & \text{if} \ \ a_i \le a_{i+1},\\
N-1 &\text{otherwise}
\end{cases}
\end{equation*}
for $1\le i\le k-1$. 
Note that if $m<0$, then $a\ensquare m \1 = D^{(-m-1)}a$.
In  this paper we extend this result to the general case.

\begin{Thm}\label{thm:basis}\sl
A linear basis of $\goth F_N(\cal B)$ is given by the set of all
monomials \fl{monom} such that 
\begin{equation}\label{fl:nojumps}
m_i - m_{i+1} \le 
\begin{cases}
\displaystyle{\sum_{j=i+1}^kN(a_i,a_j)-\sum_{j=i+2}^kN(a_{i+1},a_j)} 
& \text{if} \ \ a_i \le a_{i+1},\\[15pt] 
\displaystyle{\sum_{j=i+1}^kN(a_i,a_j)-\sum_{j=i+2}^kN(a_{i+1},a_j)-1} 
&\text{otherwise}
\end{cases}
\end{equation}
for $1\le i\le k-1$. 
\end{Thm}
Denote the set of all basic monomials by $\cal T$.
We will prove this theorem in \sec{span} and \sec{linind}.

The free vertex algebra $\goth F=\goth F_N(\cal B)$ is graded by the free
abelian semigroup $\Z_+[\cal B]$, so that the weight of the monomial
\fl{monom} is $a_1+\ldots+a_k\in \Z_+[\cal B]$. Let $\goth F_\lambda$ be
the space of all elements of $\goth F$ of weight 
$\lambda \in \Z_+[\cal B]$. 
We have the decomposition 
$\goth F = \bigoplus_{\lambda\in \Z_+[\cal B]} \goth F_\lambda$.
 There is also a grading by
$\frac 12 \Z$ such that $\deg b = -\frac 12 N(b,b)$ for 
$b\in\cal B$. Each homogeneous subspace of weight $\lambda$ is in its
turn decomposed 
$\goth F_\lambda = \bigoplus_{d\in\frac12\Z} \goth F_{\lambda,d}$,
so that  
$$
\goth F_{\lambda,d} \ensquare n \goth F_{\lambda',d'} \subset 
\goth F_{\lambda+\lambda', d+d'-n-1}
\quad \text{and}\quad
\1\in\goth F_{0,0}.
$$
Note that the parity of a homogeneous element $w\in\goth F$ is 
$p(w) \equiv 2\deg w \mod 2$ so the axiom \hyperlink{V5}{V5} holds.
The basic set $\cal T$ is homogeneous with respect to both these
gradings. Denote $\cal T_{\lambda,d} = \cal T\cap \goth F_{\lambda,d}$.

\thm{basis} 
implies the following combinatorial meaning of the dimensions of
homogeneous components of $\goth F_N(\cal B)$. Fix a weight 
$\lambda = \sum_{a\in \cal B} s_a a \in  \Z_+[\cal B]$. There is only one
basic monomial $w_{\min}(\lambda) \in \cal T$ of weight $\lambda$ 
that attains the minimal
possible degree $d_{\min}(\lambda)=\deg w_{\min}(\lambda)$. It is of the form
\begin{equation}\label{fl:wmin}
w_{\min}(\lambda)= a_1 \ensquare{m_1}( \cdots a_{k-2}  \ensquare{m_{k-2}}
(a_{k-1} \ensquare{m_{k-1}}a_k)\mbox{$\cdot\!\cdot\!\cdot$}),
\end{equation}
where $a_1\le a_2\le \ldots \le a_k$ and 
$m_i = \displaystyle{\sum_{j=i+1}^k} N(a_i,a_j)-1$ for $i\le i\le k-1$. 

Let $w\in \cal T$ be a monomial, given by \fl{monom}. Define a sequence
of integers $\eta(w) = (n_1, \ldots, n_k)\in \Z^k$ where
\begin{equation}\label{fl:eta}
n_i = \sum_{j=i+1}^k N(a_i,a_j)-1-m_j
\quad\text{for}\ \ 1\le i\le k-1, \quad
n_k = -1-m_k.
\end{equation}
By \fl{nojumps} we have $n_1\ge n_2 \ge \ldots\ge n_k$
and if $n_i=n_{i+1}$ then $a_i\le a_{i+1}$. Also,
$\sum_i n_i = \deg w - \deg_{\min}(\op{wt}w)$. We can
view $\eta(w)$ as a partition of the number $\deg w -
\deg_{\min}(\op{wt}w)$ colored by the set $\cal B$ by setting the
color of $n_i$ to be $a_i$. It is easy to see that this way we get a
one-to-one correspondence between $\cal T_{\lambda, d}$  and the set
of corresponding colored partitions of $d-d_{\min}(\lambda)$.

\begin{Cor}\label{cor:dimen}\sl
The dimension of the homogeneous component 
$\goth F_N(\cal B)_{\lambda,d}$ is equal to the number of
$\cal B$-colored partitions of $d-d_{\min}(\lambda)$ 
that contain at most
$s_a$ terms of each color $a\in\cal B$. 
\end{Cor}

In particular, these dimensions do not depend on the locality bound
$N$, up to a shift.

\section{Embedding of free vertex algebras into lattice vertex
algebras}\label{sec:embedding}
The two former sections indicate a striking similarity between lattice
vertex algebras and free vertex algebras. In this section we make this
similarity precise. 

As in  \sec{free}, let $\cal B$ be a set and let
$N:\cal B\times\cal B\to \Z$ be a symmetric function. Let
$\Lambda = \Z[\cal B]$ be the lattice generated by $\cal B$, define
the bilinear form first for $a,b\in\cal B$ by $(a|b)=-N(a,b)$
and then extend it to the whole $\Lambda$ by linearity. By the
universality property of the free vertex algebra $\goth F = \goth
F_N(\cal B)$ there is a vertex algebra homomorphism $\f:\goth F\to
V_\Lambda$ such that $\f(a)=v_a$ for each $a\in\cal B$. 
Note that $\f$ is homogeneous with respect to the double grading on $\goth F$ and $V$.
The following theorem claims that $\f$ is injective.

\begin{Thm}\label{thm:embedding}\sl
Let $\cal B\subset \Lambda$ be a linearly independent set. Then the
elements $v_a$ for $a\in\cal B$ generate a free vertex
subalgebra in $V_\Lambda$.  
\end{Thm}

We will prove this theorem in \sec{linind}. Here we deduce the
following easy fact from the mere existence of the homomorphism $\f$ and from
formula \fl{vanvb}. Let $w=w_{\min}(\lambda)\in \goth F$ be the basic minimal monomial 
of weight $\lambda = a_1 +\ldots +a_k$ given by \fl{wmin}. Then it is
easy to calculate that 
$\deg w = \deg_{\min}(\op{wt}w) = \frac 12 (\op{wt}w|\op{wt}w)$.

\begin{Prop}\label{prop:minimal}\sl
Let $w\in \goth F_\lambda$ be  monomial of weight 
$\lambda\in \Z_+[\cal B]$ such that any its submonomial $u$ has degree 
$\deg u = \deg_{\min}(\op{wt}u) = \frac 12 (\op{wt}u|\op{wt}u)$. 
Then  $\f(w)= \pm v_\lambda \neq 0$.
\end{Prop}
\begin{proof}
Use induction on the length of $w$. If $w\in \cal B$, then
$\f(w)=v_\lambda$ by the definition on $\f$. Otherwise, 
$w=u_1\ensquare m u_2$ and by the induction, $\f(u_i)=v_{\op{wt}u_i}$.
Since $\f$ is homogeneous, we have $\deg \f(w) = \deg w =
\frac12(\lambda|\lambda)$, hence $\f(w) =
v_{\op{wt}u_1+\op{wt}u_2} = v_\lambda$.
\end{proof}
Note that we do not use either of Theorems \ref{thm:basis} or
\ref{thm:embedding} in the proof of \prop{minimal}.

Clearly the  minimal monomial
$w_{\min}(\lambda)$ satisfies the assumption of this
proposition. It will follow
from Theorems \ref{thm:basis} and \ref{thm:embedding}  that 
every such monomial $w$ is \hyperlink{minimal}{minimal}  
(see the end of \sec{def}) and
is proportional to $w_{\min}(\lambda)$.

\section{Quantitative Dong's lemma}\label{sec:dong}
Let again $\cal B$ be a set with a locality function $N:\cal
B\times\cal B \to \Z$ and let $\goth F= \goth F_N(\cal B)$ be the
corresponding free vertex algebra. 
\begin{Lem}\label{lem:dong}\sl
Let $a,b,c\in\cal B$ and let $n=N(a,b)-k-1$ for some $k\ge 0$. Then 
\begin{equation}\label{fl:dong}
N(c,b\ensquare n a)=
\begin{cases}
N(a,c)+N(b,c)+k&\text{if}\quad N(b,c)>0 \ \  \text{or}\\
& \qquad 0 \le k\le -N(b,c),\\
N(a,c)&\text{if}\quad N(b,c)=0 \ \  \text{or}\\
& \qquad 0<-N(b,c)\le k.
\end{cases}
\end{equation}
\end{Lem}

If $a,b,c$ belong to some arbitrary vertex algebra $\goth A$ then
clearly the locality $N(c,b\ensquare n a)$ is bounded from above by
the right hand side of \fl{dong}. 
Note that this estimate applies also to arbitrary fields 
$a,b,c\in\F(V)$, provided they are pairwise local. Indeed,
by the \hyperlink{dong}{qualitative version} of Dong's lemma, such fields generate a
vertex algebra in $\F(V)$.

\begin{proof}
By \thm{embedding} it is enough to prove the lemma for the case
when $a = v_\alpha$, $b=v_\beta$ and $c=v_\gamma$ are vacuum
elements in a lattice vertex algebra $V_\Lambda$ for some vectors  
$\alpha,\beta,\gamma\in\Lambda$.  

Let $m = -(\alpha+\beta|\gamma)+k-j-1$ for some $j \ge 0$.
Using \fl{vanvb} together with the formulas
\begin{gather*}
v_\beta \ensquare n \big(Du\big) = 
n\, v_\beta \ensquare{n-1} u + 
D\big( v_\beta \ensquare n u\big),\\
v_\beta \ensquare n \big(\alpha(-k)\,u\big) = 
-(\alpha|\beta)\, v_\beta \ensquare{n-k} u + 
\alpha(-k)\,\big( v_\beta \ensquare n u\big),
\end{gather*}
for $\alpha\in\goth h$,\  $\beta \in\Lambda$ and $u\in V_\Lambda$, we obtain 
\begin{multline}
c \ensquare{m} (b \ensquare n a) = 
\e(\alpha,\beta)\e(\alpha,\gamma)\e(\beta,\gamma)
\sum_{i=k-j}^k\binom{k-(\beta|\gamma)-j-1}i\\
\times\big(D-\alpha(-1)\big)^{(k-i)}
\big(D-(\alpha+\beta)(-1)\big)^{(j-k+i)}\,  
v_{\alpha+\beta+\gamma}.\notag
\end{multline}
Let $j_{\min}$ be the minimal value of $j$ such that 
$c \ensquare{m} (b \ensquare{n} a)\neq 0$. We have 
$$
N(c, b\ensquare n a) = -(\alpha+\beta|\gamma)
+k-j_{\min}=N(a,c)+N(b,c)+k-j_{\min}.
$$
It follows that if $(\beta|\gamma)<0$ or if 
$(\beta|\gamma)\ge k\ge 0$, then   $j_{\min}=0$; 
if $(\beta|\gamma)=0$, then $j_{\min} = k$; finally, 
if $k>(\beta|\gamma)>0$, then $j_{\min}=k-(\beta|\gamma)$
and the statement follows.
\end{proof}

For monomials of length more than three the analogous estimate is more
subtle. Note however, that by \cor{dimen}, if 
$w\in \goth F$ is a vertex monomial in $\cal B$ such that 
$\deg w < \deg_{\min}(\op{wt}w)$ then $w=0$. So we get:
\begin{Prop}\label{prop:dong}\sl
Let $w=a_1\ensquare{n_1}\cdots\ensquare{n_{l-1}}a_l \in \goth F$, \ 
$a_i\in \cal B, \ n_i\in \Z$, be a vertex monomial. If
\begin{equation}\label{fl:ineq}
\sum_{i=1}^{l-1} n_i > \sum_{1\le i<j\le l}N(a_i,a_j)-l+1,
\end{equation}
then $w=0$.
\end{Prop}
\prop{minimal} shows that sometimes the estimate \fl{ineq} is the best
possible. For $l=3$ \lem{dong} gives in general a stronger estimate.

\section{Locality function}\label{sec:locfun}
Let $\goth L$ be a \hyperlink{conformal}{conformal algebra} (see
\sec{conformal}) generated by a finite set 
$\cal B\subset \goth L$. In \cite{universal} the following integer
function was defined. Consider all monomials 
$w=a_1\ensquare{n_1}\cdots \ensquare{n_{l-1}}a_l\in\goth L$, \ 
$a_i\in\cal B$, \ $n_i \in \Z_+$.  Let $S(l)$ be the maximal 
possible value of $\sum_{i=1}^{l-1} n_i$ such that $w\neq0$. We call $S$
the \hypertarget{locfun}{locality function} of $\goth L$. It depends
on the generating set 
$\cal B$, however, if we take another generating set, then the growth
of $S$ can change  at most by a linear term. It was shown in
\cite{universal} that if $\goth L$ is embeddable into an {\it
associative conformal algebra}, then $S(l)$ must have at most linear
growth. 

Clearly, similar function can be defined for a finitely generated
vertex algebra.

Now let $\cal B$ be a finite set with a locality bound 
$N:\cal B\times\cal B\to
\Z_+$ that takes only non-negative values. Let $\goth F = \goth
F_N(\cal B)$ be the corresponding free vertex algebra (see
\sec{free}). Let $\goth L \subset \goth F$ be the conformal algebra
generated by $\cal B$. It follows from the results of \cite{freecv}
that $\goth L$ is a free conformal algebra. It was shown in
\cite{universal} that if $N\not\equiv0$ then the
locality function $S(l)$ of $\goth L$ has at least
quadratic growth, and hence $\goth L$ cannot be embedded into an associative
conformal algebra. It was also conjectured that $S(l)$ has exactly
quadratic growth. Combining Propositions \ref{prop:minimal} and 
\ref{prop:dong}, we see that this conjecture is indeed true.

\begin{Cor}\label{cor:locfun}\sl
The locality function of the free conformal algebra generated by a
set $\cal B$ with locality bound $N$ is  
$$
S(l) = \frac{l(l-1)}2\max_{a,b\in\cal B}N(a,b) -l+1.
$$
\end{Cor}

\section{Actions of conformal algebras on free vertex
algebras}\label{sec:confact}
Return to the setup of \sec{embedding}. 
Let $\goth F = \goth F_N(\cal B)$ be the free vertex algebra generated
by a set $\cal B$ with a locality bound $N:\cal B\times\cal B\to
\Z$. Let $\Lambda = \Z[\cal B]$ and let $\f:\goth F\to
V_\Lambda$ be the injective homomorphism provided by \thm{embedding}, 
such that $\f(a)=v_a$ for $a\in\cal B$.

We are interested in \hyperlink{confder}{conformal derivations} 
$\alpha \in cder\,\goth F$,
which are homogeneous of weight 0 and of some degree $m+1\in \Z$ (see
\sec{conformal}). This 
means that for each $b\in \cal B$ we have 
$\alpha(n)b = f_n(b)D^{m-n} b$ 
for some functions $f_n:\cal B \to \Bbbk$ when $0\le n \le m$, 
and $\alpha(n)b = 0$ when $n>m$. Recall that a conformal derivation is
uniquely defined by its action on the generators. 

We study the case when $m=0$ in a greater detail, because it will be used
later in \sec{linind}.
Recall from \sec{lattice} that $V_0$ contains the conformal Heisenberg
algebra $\goth H$. This algebra
acts on $V_\Lambda$ by commuting inner conformal derivations $\tl h_+$ for 
$h\in\goth h = \Lambda\otimes \Bbbk$ such that $N(\tl h_+, \tl g_+)=0$.

\begin{Lem}\label{lem:H+}\sl
\begin{enumerate}
\item\label{H:f}
Let $f:\cal B \to \Bbbk$ be an arbitrary function.
Then there is a unique conformal derivation $\alpha_f\in
cder\,\goth F$ such that 
$N(\alpha_f, \cal B) = 1$
and $\alpha_f(0)b = f(b)\,b$ for all $b\in \cal B$. 
For another function $g:N\to \Bbbk$ we
have $N(\alpha_f, \alpha_g)=0$.
\item\label{H:h}
Let $h\in\goth h$ is such that $f(b) = (h|b)$ for all $b\in\cal B$.
Then $\f\big(\alpha_f(n)x\big) = h(n)\f(x)$ for every $x\in\goth F$ and $n\ge0$.
\end{enumerate}
\end{Lem}

\begin{proof}\ 
\ref{H:f}.\
Consider the derivation $\alpha_f(n)$,\ $n\ge 0$, of the associative
algebra $\ol{\Bbbk\<\cal X\>}$ which acts on a generator $b(m)\in \cal X$ 
by $\alpha_f(n)\big(b(m)\big) = f(b)\, b(m+n)$, see \sec{free} for the
notations. It is easy to see that 
$\alpha_f(n)$ preserves the ideal $I\subset \ol{\Bbbk\<\cal X\>}$ generated
by the locality relations \fl{locoef},  hence $\alpha_f(n)$
is a derivation of the algebra $U=  \ol{\Bbbk\<\cal X\>}/I$. Also, we have 
$\ad D{\alpha_f(n)}=-n\,\alpha_f(n-1)$ where $D$ is the derivation
of $U$ defined by $Db(m) = -m\,b(m-1)$. Since $\alpha_f(n)$ can only
increase the indices, it preserves the left ideal $U\cal X_+\subset U$,
therefore $\alpha_f(n)$ acts on the free vertex algebra 
$\goth F = U/U\cal X_+$.

\smallskip\noindent
\ref{H:h}.\ 
Both $\alpha_f$ and $\tl h_+$
are conformal derivations, therefore they are uniquely determined by
their action on the generators $\cal B$. But for every $b\in \cal B$ 
we have $\alpha_f(n)b = \delta_{n,0}\, (h|b)\,b$ and 
$h(n)v_b = \delta_{n,0}\, (h|b)\,v_b$, and the claim follows.
\end{proof}

Let us now consider the case when $m=1$. A derivation of weight 0 and
degree 2 is determined
by a pair of functions $f_0, f_1:\cal B\to \Bbbk$. One can show that
in order for this derivation to be well defined $f_0$ must be
constant. Take $f_0 \equiv 1$ and an arbitrary
function $f=f_1:\cal B\to \Bbbk$. This defines a conformal derivation 
$\omega_f\in cder\,\goth F$ such that $N(\omega_f, \cal B) = 2$
and $\omega_f(0)b = Db$, \   $\omega_f(1)b = f(b)\,b$ for every 
$b\in \cal B$. For another function $g:\cal B\to \Bbbk$ we
have $N(\alpha_f, \omega_g)=N(\omega_f, \omega_g)=2$ and 
\begin{equation*}
\omega_f\ensquare 0 \alpha_g = \frac d{dz}\alpha_g, \ \
\omega_f\ensquare 1 \alpha_g = \alpha_g, \ \
\omega_f\ensquare 0 \omega_g =  \frac d{dz}\omega_g, \ \
\omega_f\ensquare 1 \omega_g =  2\omega_g. 
\end{equation*}
In particular each $\omega_f$ generates an action of the Virasoro conformal
algebra $\goth{Vir}$. Note that the products 
$\alpha_f\ensquare n \omega_g$ can be calculated using \fl{qs}.

All these statements are proved in a way similar to the proof of \lem{H+}.
Using \thm{present} in \sec{present} one can show that the derivation $\omega_f$ can be
extended from $\f(\goth F)$  to $V_\Lambda$ if and only
if $f(b) = \frac 12(b|b)$. If $V_\Lambda$
has the Virasoro element $\omega \in V_0$, then this extension
coincides with the inner conformal derivation $\omega_+$ defined by $\omega$.

\begin{Rem}
This shows that though free vertex algebras are not vertex operator
algebras, they have an action of the conformal Virasoro algebra
compatible with the $\frac12\Z$-grading. 
\end{Rem}

However, derivations of higher degree  exist only in some very special cases.

\subsubsection*{Bozon--fermion correspondence}
Suppose that $\cal B = \{a\}$ consists of only one element with 
$N(a,a)=-1$. Then $\Lambda = \Z$ and the form is given by
$(m|n)=mn$. The corresponding 
lattice vertex algebra $V_\Z$ is the main object of the so-called 
bozon--fermion correspondence, see e.g. \cite{kac2}.
It is well known that the Heisenberg vertex algebra 
$V_0 \subset V_\Z$ contains a
conformal algebra $\^{\goth W}$ spanned over $\Bbbk[D]$ by the
elements $p_m = v_{-1}\ensquare{-m-1}v_1$, \ $m\in \Z_+$. The
multiplication table is  
\begin{multline*}
p_m\ensquare k p_n = 
\binom{m+n-k}m\, p_{m+n-k} \\
-\sum_{s=0}^{m-k} (-1)^{k+s}  \binom{m+n-k-s}n\, D^{(s)} p_{m+n-k-s}
+\delta_{k,m+n+1}(-1)^m \1.
\end{multline*}
In particular, $p_0 = -\tl a$ generates the Heisenberg conformal
algebra $\goth H\subset \^{\goth W}$, and $p_1 = \frac12 \tl a \ensquare{-1} \tl a
- \frac 12 D\tl a$ generates the Virasoro conformal algebra 
$\goth{Vir}\subset \^{\goth W}$.

The \hyperlink{coeff}{coefficient algebra} of $\^{\goth W}$ is a central extension
$\^W=W\oplus \Bbbk\c$ 
of the Lie algebra $W=\Bbbk\<p,t^{\pm1}\,\big|\,\ad tp=1\,\>\lie$
of differential operators on the circle. The subalgebra
$W_+\subset \^W$  spanned by 
the coefficients $p_m(n)$ for $n\ge 0$ can be identified with
the Lie algebra $\Bbbk\<p,t\,\big|\,\ad tp=1\,\>\lie$  of
differential operators on the disk so that $p_m(n) =
\frac1{m!}p^mt^n$, see \cite{fkrw,bozfer}.

Using \fl{vanvb} and the fact that $\e=1$ we calculate that 
\begin{equation}\label{fl:p}
p_m(n)v_1 = \begin{cases}
(-1)^m\, D^{(m-n)} v_1&\text{if}\ \ m\ge n,\\
0&\text{otherwise}.
\end{cases}
\end{equation}
Therefore, the inner conformal derivations $(p_m)_+\in cder\,V_\Z$ 
are homogeneous of weight 0 and degree $m+1$.
We will now show that $(p_m)_+$ preserves $\f(\goth F)$, so by
\thm{embedding} it acts on $\goth F$.

Each vacuum vector $v_k\in V_\Z$ is a lowest weight vector for the Lie
algebra $W_+$, meaning that for a homogeneous $g\in W_+$ one has 
$gv_k=0$ if $\deg g >0$ and $gv_k = \lambda(g)v_k$ for some weight
$\lambda\in (W_+)_0^*$ if $\deg g=0$. It is
shown in \cite{bozfer} that $v_k$ generates an irreducible lowest
weight $W_+$-module $U(W_+)v_k\subset V_k$ of the lowest weight $\lambda_k$ defined
by $\lambda_k\big(p_m(m)\big) = (-1)^m \binom{m+k}k$. It is also proved in
\cite{bozfer} that dimension of
the homogeneous component of $U(W_+)v_k$ of degree $d+k^2/2$ is equal to
the number of partitions of $d$ into at most $k$ terms, i.e. is
exactly $\dim \goth F_{k,d+k^2/2}$. This motivates the following
theorem.

\begin{Thm}\label{thm:bozfer}\sl
$\f(\goth F_k) = U(W_+)v_k$ for every $k\ge 1$.
\end{Thm}

\begin{proof} \ 
Using the dimension argument it is enough to prove only that 
$U(W_+)v_k\subseteq\f(\goth F_k)$, though in fact the proof of the other
inclusion is no more difficult. 

Since $v_k \in \f(\goth F)$, it is enough to show that 
$W_+\f(\goth F)\subset \f(\goth F)$. The Lie algebra $W_+$ is spanned by $p_m(n)$ 
for $m,n\ge 0$.  From  \fl{p} it follows that $p_m(n)v_1 \in
\f(\goth F)$. Since $v_1$ is a generator of $\f(\goth F)$, it follows
that  each $p_m(n)$ preserves $\f(\goth F)$. 
\end{proof}

\begin{Rem}
One can show that homogeneous conformal derivations of weight 0 and
degree $m>2$ of a free vertex algebra $\goth F$ exist if and only if $\goth
F$ is a tensor product of free vertex algebras, generated by a single
generator $a$ such that either $N(a,a)=0$ or $N(a,a)=-1$.
\end{Rem}
\section{Presentation of lattice vertex algebras in terms of generators
and relations}\label{sec:present}
In this section we start with a lattice vertex algebra $V_\Lambda$
corresponding to an integer lattice $\Lambda$. 
Choose a $\Z$-basis $\Pi$  of $\Lambda$.
We may assume that $\e(a,a)=1$ for every $a\in \Pi$.
By \thm{embedding} the vertex subalgebra of $V_\Lambda$ generated by 
$\set{v_a}{a\in\Pi}$ is isomorphic to a free vertex algebra. 
On the other hand, it follows from \fl{vanvb} that if we complete
these generators by $v_{-a}$'s, then the resulting vertex algebra will
be the whole $V_\Lambda$.

\begin{Thm}\label{thm:present}\sl 
The lattice vertex algebra $V_\Lambda$ is presented by generators 
$\set{v_a}{a\in\pm\Pi}$ with locality bound given by
$N(v_a,v_b)=-(a|b)$ and relations 
$\set{v_a \ensquare{(a|a)-1}v_{-a}=\1}{a\in\Pi}$.
\end{Thm}
\noindent
Note that the quasisymmetry identity \fl{qs} implies that also 
$v_{-a} \ensquare{(a|a)-1}v_{a}=\1$ for $a\in\Pi$.

For the proof we need the following lemma.

\begin{Lem}\label{lem:present}\sl
The lattice vertex algebra $V_\Lambda$ is
presented by generators $\tl  a$ and  $v_{\pm a}$ for 
$a\in \Pi$ with locality
bound given by
\begin{equation}\label{fl:presloc}
N(v_a, v_b) = -(a|b),\quad
N(\tl  a, \tl  b) = 2,\quad
N(\tl  a,v_b)=1,
\end{equation}
and the following set of relations: 
\begin{itemize}
\item[(i)]
$\tl  a \ensquare 0 \tl  b = 0, \ \tl  a \ensquare 1 \tl  b =
(a|b)\,\1$,
\item[(ii)]
$\tl  a\ensquare 0 v_b = (a|b)\,v_b$,
\item[(iii)]
$v_a \ensquare{(a|a)-1}v_{-a} = \1$,
\item[(iv)]
$Dv_a = \tl  a\ensquare{-1}v_a$,
\end{itemize}
for 
any $a,b \in \pm\Pi$.
\end{Lem}
The relations (i) mean that 
$\tl  a$'s generate the Heisenberg vertex algebra $\goth H
\subset V_\Lambda$. The relations (ii) say that each $v_b$ is a
vacuum vector for the representation of the corresponding 
Heisenberg Lie algebra of weight given by 
$a \mapsto (a|b)$. 

We first show how this lemma implies the theorem.


\begin{proof}[Proof of \thm{present}]
Consider a set of generators
$\cal B =\set{X_a, X_{-a}}{a\in\Pi}$ with the  
locality bound $N:\cal B\times\cal B\to\Z$ given by 
$N(X_a,X_b)=-(a|b)$. 
Let $\goth F = \goth F_N(\cal B)$ be the corresponding free vertex
algebra. For $a\in \pm\Pi$ denote $K_a = X_a \ensquare{(a|a)-1}X_{-a}$,\,
$H_a = X_a \ensquare{(a|a)-2}X_{-a}\in\goth F$.
Let $\goth I \subset \goth F$ be the ideal generated by the
relations $K_a=\1$ for 
$a\in \Pi$. We have to show that the quotient vertex algebra 
$\goth A = \goth F/\goth I$ is isomorphic to $V_\Lambda$ so that 
$X_a \mapsto v_a$ under the projection map
$\pi:\goth F \to \goth F/\goth I = \goth A$.

By \lem{present} it
is enough to show that the elements 
$\tl a = \pi(H_a)$,\, $v_a = \pi(X_a) \in \goth A$ satisfy
the identities (i)--(iv) above. Note that the locality of these
elements is already as prescribed by \fl{presloc}, and (iii) holds by the
assumption. The identities (i), (ii) and (iv) follow from  the following
identities that hold in the free vertex algebra $\goth F$ for all 
$a,b\in\pm\Pi$:
\begin{equation}\label{fl:idF}
\begin{split}
H_a\ensquare k H_b &= (a|b)\, K_a \ensquare{k-2} K_b,
\quad\text{for}\ \ k=0,1,\\
H_a \ensquare 0 X_b &= (a|b)\, K_a\ensquare{-1}X_b,\\
H_a\ensquare{-1}X_a &= X_a \ensquare{-2}K_a + (a|a)\, K_a \ensquare{-2}X_a.
\end{split}
\end{equation}

It is quite easy to prove these identities by straightforward
calculations using only the construction of $\goth F$ given in
\sec{free}. However, we can simplify these calculations even further
by using \thm{embedding}.  
Consider the set $\bar\Pi=\set{a,\bar a}{a\in\Pi}$ and extend the symmetric 
bilinear form to $\bar\Lambda=\Z\left[\bar\Pi\right]$ by 
$(a|\bar b) = - (a|b)$,\  
$(\bar a|\bar b) = (a|b)$. In the same way we
extend the cocycle $\e$, assuming that $\e(a,\bar a)=\e(\bar a, \bar a)=1$ for
$a\in\Pi$. Let  $V_{\bar\Lambda}$ be the
corresponding lattice vertex algebra and let 
$\f:\goth F\to V_{\bar\Lambda}$ be the injective homomorphism constructed in
\sec{embedding}, such that $\f(X_a) = v_a$, \ 
$\f(X_{-a}) = v_{\bar a}$ for 
$a\in \Pi$. Let us, for example, prove the first identity in \fl{idF},
the other two are proved in the same way.

Take some $a,b \in \Pi$. By \fl{vanvb} we have 
$\f(K_a) = v_{a+\bar a}$ and 
$\f(H_a) = a(-1)v_{a+\bar a}$, so, using \fl{assoc} and \fl{comm}, 
we get for any $k\in \Z_+$:
\begin{align*}
\f(H_a)\ensquare k \f(H_b) 
&= \big(a(-1)v_{a+\bar a}\big) \ensquare k 
\big(b(-1)v_{b+\bar b}\big) \\
&=\sum_{s\ge 0} a(-s-1)v_{a+\bar a}(k+s)b(-1)v_{b+\bar b}\\
&+\sum_{s\ge0}v_{a+\bar a}(k-s-1)a(s)b(-1)v_{b+\bar b}.
\end{align*}
Since $(b|a+\bar a)=(a|b+\bar b)=0$ we have 
$\ad{b(i)}{v_{a+\bar a}(j)}=
a(0)v_{b+\bar b}=0$, so the only non-zero term in the above comes
from the second sum when $s=1$ and we precisely get the first identity
in \fl{idF}.
\end{proof}

\subsubsection*{Proof of \lem{present}}
Let $\goth A$ be the vertex algebra presented by the generators 
$\set{\tl  a, v_{\pm a}}{a\in\Pi}$ with locality
bound \fl{presloc} and relations (i)--(iv). Since these relations hold
in the lattice vertex algebra $V_\Lambda$, there is a projection 
$\pi:\goth A\to V_\Lambda$. 

Let $w = v_{b_1}\ensquare{n_1}\cdots
\ensquare{n_{k-1}}v_{b_k}\in \goth A$, \ 
$b_i\in \pm \Pi$, be a vertex monomial. 
Denote by $\Xi(w)$ the set of all vertex monomials involving
the same indeterminates $v_{b_i}$'s in the same order as in $w$ and
with the same arrangement of parentheses. Let $\mu(w)\in\Xi(w)$ be the
unique \hyperlink{minimal}{minimal} monomial in $\Xi(w)$, see
\sec{def}. As in \prop{minimal}, we have  
$\pi\big(\mu(w)\big) = \pm v_{\op{wt}w}\in V_\Lambda$. 

Recall from \sec{lattice} that $V_\Lambda$ is a module over the 
\hyperlink{heisen}{Heisenberg algebra} $H = \goth h\otimes
\Bbbk[t,t\inv]\oplus \Bbbk\c$ for $\goth h = \Lambda \otimes_\Z
\Bbbk$. As before we denote $a\otimes t^n=a(n)$. 
The algebra $H = H_-\oplus H_+$ is decomposed into a direct
sum of subalgebras $H_-= \set{h(n)}{h\in\goth h,\ n<0}\oplus \Bbbk\c$
and $H_+= \set{h(n)}{h\in\goth h,\ n\ge0}$.
Let $\^H_- = H_- \oplus \Bbbk D$ be the 
Lie algebra obtained from $H_-$ by adjoining the derivation 
$D$ acting by $h(n)\mapsto -n\,h(n-1)$ and $D\c=0$. 
Let $\^H = \^H_- \oplus H_+$.
Then $V_\Lambda$ is a module over $\^H$.

For a linear combination $h = \sum_{a\in\Pi} s_a a\in \goth h$, \
$s_a\in\Bbbk$, denote $\tl h = \sum_{a\in\Pi} s_a \tl  a\in \goth A$. 
We observe that $\goth A$ is also a module over $\^H$ with the action
given by $h(n)x=\tl h \ensquare n x$, and the vertex algebra
homomorphism $\pi:\goth A\to V_\Lambda$ is also a $\^H$-module
homomorphism. 

\begin{Lem}\label{lem:maximal}\sl
Let $w=\mu(w)\in\goth A$ be a right-normed minimal monomial of weight 
$\lambda \in \Lambda$. Then 
\begin{enumerate}
\item\label{max:N}
For any $a \in \pm\Pi$, we
have $N(v_a, w) = -(a|\lambda)$. 
\item\label{max:D}
$Dw=\lambda(-1)w$.
\end{enumerate}
\end{Lem}
\begin{proof}
We prove this lemma by induction on the length $l$ of $w$. If $l=1$
then $w = v_\lambda$ and statements \ref{max:N} and \ref{max:D}
are true by assumption.

Assume $l>1$. Then $w = v_b \ensquare n u$ for some other maximal 
right-normed monomial $u$ of weight $\lambda'= \lambda-b$. 
By the induction hypothesis, 
$N(v_b, u)= -(b, \lambda') = n+1$. Using \lem{dong}, 
we have: $N(v_a, v_b \ensquare n u) =
-(a|b)-(a|\lambda')-(b|\lambda')-n-1=
-(a|\lambda)$. This proves \ref{max:N}. 

In order to prove \ref{max:D} we first calculate, using \fl{assoc}:
\begin{align*}
\big(b(-1)v_b\big)\ensquare n u &=
\sum_{s<0}b(s)\big(v_b\ensquare{n-s-1}u\big) +  
\sum_{s\ge0}v_b\ensquare{n-s-1}\big(b(s)u\big)\\
&= b(-1)\big(v_b\ensquare n u\big) + 
(b|\lambda')\, v_b\ensquare{n-1}u.
\end{align*}
Now we have, using the induction hypethesis and \fl{comm},
\begin{align*}
D\big(v_b\ensquare n u\big) &= 
\big(b(-1)v_b\big) \ensquare n u + 
v_b \ensquare n \big(\lambda'(-1)u\big)\\
&= b(-1)\big(v_b\ensquare n u\big) + 
(b|\lambda')\, v_b\ensquare{n-1}u \\
&\qquad - (b|\lambda')\, v_b\ensquare{n-1}u + 
\lambda'(-1)\big(v_b\ensquare n u\big) \\
&= \lambda(-1)\big(v_b\ensquare n u\big).
\end{align*}
\end{proof}

An element $x\in\goth A$ is called vacuum, if for every
$h\in\goth h$ one has $h(n)x=0$ if $n>0$, and $h(0)x = \lambda(h)x$
for some $\lambda\in\goth h^*$. It is easy to see that a minimal
monomial $w=\mu(w)\in\goth A$ is a vacuum element. 
We prove next that the minimal monomials generate
$\goth A$ over the extended Heisenberg algebra $\^H$. It will follow
then by the representation theory of Heisenberg algebras (see
e.g. \cite{kac1}) that as a module over $\^H$ the vertex algebra 
$\goth A = \bigoplus_{w=\mu(w)}U(\^H)w$ is decomposed into a direct
sum of irreducible highest weight modules generated by minimal monomials.
Recall that by \fl{assoc} a vertex algebra is spanned by right-normed
monomials in its generators.

\begin{Lem}\label{lem:vacuum}\sl
Let $w\in\goth A$ be a right-normed monomial in the variables\break
$\set{v_a}{a\in \pm \Pi}$. Then 
$w\in U(\^H_-)\mu(w)$.
\end{Lem}
\begin{proof}
As before we prove this lemma by induction on the length $l$ of $w$. If
$l=1$, then $w=\mu(w) = v_a$. Otherwise,
$w = v_a \ensquare n u$ for some right-normed monomial $u$. By
induction, we have $u = g\, \mu(u)$ for some $g \in
U(\^H_-)$. Applying the formulas
\begin{gather*}
v_a \ensquare n \big(Du\big) = 
n\, v_a \ensquare{n-1} u + 
D\big( v_a \ensquare n u\big),\\
v_a \ensquare n \big(h(-k)u\big) = 
-(a|h)\, v_a \ensquare{n-k} u + 
h(-k)\big( v_a \ensquare n u\big),\quad h\in\goth h,
\end{gather*}
which follow from 
\hyperlink{V3}{V3} and \fl{comm},
we can express $w$ as 
$$
w=\sum_{j\in \Z} g_j\, v_a\ensquare j
\mu(u),\qquad g_j \in U(\^H_-).
$$ 
Thus we have reduced the lemma
to the case when $w=v_a\ensquare n u$ for a minimal monomial 
$u=\mu(u)$. 

Let $\lambda = \op{wt}u\in \Lambda$. By \lem{maximal}\ref{max:N}, we
have $N(v_a, u) = -(a|\lambda)$. Assume that $w$ is not
maximal, then $n =-(a|\lambda) - k-1$ for some $k>0$. Using  
\lem{maximal}\ref{max:D} and the formulas above, we get 
\begin{align*}
\lambda(-1)\big(v_a \ensquare{n+1}u\big) &=
(a|\lambda)\,v_a\ensquare n u
+ v_a \ensquare{n+1}Du \\
&= (a|\lambda)\,v_a\ensquare n u
+(n+1)\,v_a\ensquare n u 
+D\big(v_a \ensquare{n+1}u\big)\\
&= D\big(v_a \ensquare{n+1}u\big) 
- k \,v_a\ensquare n u.
\end{align*}
Therefore, $v_a \ensquare n u = \frac 1k \big(D-\lambda(-1)\big) 
\big(v_a \ensquare{n+1}u\big)$, and the lemma follows.
\end{proof}

Note that here we already recover the formula \fl{vanvb}, up to a
sign.

\smallskip

As in the proof of \thm{present}, consider the free vertex algebra
$\goth F=\goth F_N(\cal B)$ generated by the set 
$\cal B =\set{X_a}{a\in\pm\Pi}$. Recall that $\goth F$ is graded by
the semilattice  $\bar \Lambda_+ = \Z_+\left[\bar\Pi\right]$, where 
$\bar\Pi=\set{a,\bar a}{a\in\Pi}$, so that 
$\op{wt}X_a = a$, \ $\op{wt}X_{-a} = \bar a$ for 
$a \in \Pi$. Let $\psi:\goth F\to \goth A$ be the vertex algebra
homomorphism such that $\psi(X_a) = v_a$ for $a\in \pm\Pi$.

By \cor{dimen}, all monomials in $\goth F$ of a weight 
$\lambda\in\bar \Lambda_+$ having the
minimal possible degree $d_{\min}(\lambda)$ are proportional to
$w_{\min}(\lambda)$, given by \fl{wmin}. If $\lambda$ contains a pair 
$a, \bar a$, then by permuting the variables in $w_{\min}(\lambda)$ we
obtain a minimal monomial of the form
$u\ensquare n \big(v_a \ensquare{(a|a)-1}v_{-a}\big)\in\goth F_\lambda$.
Applying (iii) we get that 
$\psi\big(w_{\min}(\lambda)\big)$ is proportional to 
$\psi\big(w_{\min}(\lambda-a-\bar a)\big)$ in $\goth A$.
It follows that  minimal monomials in $\goth A$ in the variables 
$\set{v_a}{a\in \pm \Pi}$ are
parametrized by the lattice $\Lambda$. Combining this with
\lem{vacuum}, we see that 
the vertex algebra $\goth A$ is decomposed into a direct sum $\goth A =
\bigoplus_{\lambda \in \Lambda}\goth A_\lambda$ of irreducible highest
weight $\^H$-modules $\goth A_\lambda$. 
Since $V_\Lambda$ has the same decomposition, the projection 
$\pi:\goth A\to V_\Lambda$ must be an isomorphism.
This finishes the proof of \lem{present}.
 
\begin{Rem}
In fact it is not very difficult to make the last argument without a
reference to \cor{dimen}, thus rendering the whole proof of \thm{present}
more or less independent on Theorems \ref{thm:basis} and \ref{thm:embedding}.
\end{Rem}

\begin{Rem}
One can prove \lem{present} in a different way, using the construction
of $V_\Lambda$ by vertex operators, and the so-called method of
$Z$-algebras, that originates to the work of Lepowsky and Wilson
\cite{lepwil81}, see also  \cite{dong,flm,lixu}. 
This remark is due to C.~Dong.
\end{Rem}

\section{Proof of theorem \lowercase{\ref{thm:basis}}, part I}\label{sec:span}
Let as before $\cal B$ be a set with a symmetric locality bound 
$N:\cal B\times \cal B\to \Z$ and let $\goth F=\goth F_N(\cal B)$ be
the corresponding free vertex algebra. 
Let $\cal X = \set{a(n)}{a\in\cal B,\ n\in\Z}$. As in \sec{free}, consider
the completion $\ol{\Bbbk\<\cal X\>}$ of the free associative algebra
generated by $\cal X$. Every element $g\in \ol{\Bbbk\<\cal X\>}$ is 
a linear combination, possibly infinite, of words $\cal X^*$ in the alphabet $\cal X$.

 By \sec{free} there is a map $\rho:\ol{\Bbbk\<\cal X\>}\to
\goth F$, defined by $a_1(n_1)\cdots a_k(n_k)\mapsto 
a_1 \ensquare{n_1} (a_2  \ensquare{n_2}
\cdots (a_k \ensquare{n_k}\1)\mbox{$\cdot\!\cdot\!\cdot$})$.
Under this map the set $\cal T$, described by \thm{basis}, can be
identified with the set of 
words $a_1(n_1)\cdots a_k(n_k)\in \cal X^*$ satisfying the condition
\fl{nojumps}.
In this section we prove that $\rho\cal T$ spans  $\goth F$ over $\Bbbk$.

Let $w\in \cal X^*$  and 
$f\in \ol{\Bbbk\<\cal X\>}$. A  rule
$w\xra{\ r\ }f$ on $\ol{\Bbbk\<\cal X\>}$ is a partially
defined linear map $\ol{\Bbbk\<\cal X\>}\to\ol{\Bbbk\<\cal X\>}$
which is applicable to an element 
$g\in \ol{\Bbbk\<\cal X\>}$ if $w$ occures in the decomposition of $g$ into
a linear combination of words from $\cal X^*$. The result $h$ of
the application of $r$ to $g$ is obtained by substituting $f$
instead of $w$; in this case we write  $g\xra{\ r\ }h$. For a set of rules 
$\cal R$ denote by $M(\cal R) = \set{w\in \cal X^*}{\exists\, w\ra f\in
\cal R}$ the set of words to which at least one rule from $\cal R$ is
applicable.

Now we construct a set of rules $\cal R$ that will correspond to the locality
relations \fl{locoef} in the generators $\cal B$. Let 
$w=a_1(n_1)\cdots a_k(n_k)\in \cal X^*$. Define 
\begin{equation}\label{fl:m}
\begin{split}
m_j &= \sum_{i=j+1}^k N(a_j,a_i) - \sum_{i=j+2}^k N(a_{j+1},a_i)\quad
\text{for}\ 1\le j\le k-2,\quad\\
m_{k-1}&=N(a_{k-1},a_k).
\end{split}
\end{equation}
Suppose there is $1\le j\le k-1$ such that $n_j - n_{j+1} > m_j$ or 
$n_j - n_{j+1} = m_j$ and $a_j>a_{j+1}$. Then $\cal R$ contains the
rule 
\begin{equation*}
w\ra a_1(n_1)\cdots a_{j-1}(n_{j-1})\,f_j\, 
a_{j+2}(n_{j+2})\cdots a_k(n_k)
\end{equation*}
for
\begin{multline*}
f_j =-\sum_{s\ge1}(-1)^s\binom Ns\,
a_j(n_j-s)a_{j+1}(n_{j+1}+s)\\
+(-1)^{p(a_j)p(a_{j+1})} 
\sum_{s\le N} (-1)^s \binom N{N-s}\,
a_{j+1}(n_{j+1}+s)a_j(n_j-s),
\end{multline*}
where $N=N(a_j,a_{j+1})$.

We need another set $\cal Q$ of rules on  $\ol{\Bbbk\<\cal X\>}$. 
Consider a sequence of elements $a_1, \ldots, a_k \in \goth F$. By the
\hyperlink{dong}{qualitative Dong's lemma} there is $S = S(a_1, \ldots, a_k)\in\Z$ such
that any vertex monomial $w =
a_1\ensquare{n_1}\cdots\ensquare{n_{k-1}}a_k\in\goth F_N(\cal B)$ is equal
to 0 whenever $\sum_{i=1}^{k-1}n_i \ge S$.

\begin{Rem}
It will follow that 
$S(a_1, \ldots, a_k) = \sum_{1\le i<j\le k}N(a_i,a_j)-k+2$, 
see \prop{dong}. We could prove this right now, but
our argument does not make use of this sharp estimate.
\end{Rem}

Define  $\cal Q$ to be the set of all rules $w\ra 0$ such that 
$w=a_1(n_1)\cdots a_k(n_k)$ has the following property: there is 
$1\le j\le k$ such that 
$\sum_{i=j}^k n_i \ge S(a_j, \ldots, a_k, \1)$. In particular, if 
$n_k\ge 0$, then $w\ra 0\in\cal Q$. 
We observe that the set $\cal T$ is exactly the set of all terminal
monomials with respect to the rules $\^{\cal R}=\cal R\cup \cal Q$, i.e. the
set of monomials to which these rules can not be applied. On the other hand,
for any rule $w\ra f\in \^{\cal R}$ the identity $\rho(w-f)=0$ holds 
in $\goth F$. Therefore, all we have to prove is that every word 
$w\in \cal X^*$ can be 
reduced to a linear combination of terminal words by a (possibly,
infinite) number of applications of the rules $\^{\cal R}$. Thi
sterminal linear combination is necessary finite, because the rules
$\^{\cal R}$ preserve weight and 
degree, and there are only finitely many
monomials in $\cal T$ of given weight and degree, see \cor{dimen}.

We will apply the rules  $\^{\cal R}$  to a word 
$w\in \cal X^*$ is such a way that $\cal Q$ has priority over $\cal R$,
in other words, we apply $\cal R$ only if $\cal Q$ can no longer be
applied. Denote by $\Omega(w)$ the set of all words 
$u \not\in M(\cal Q)$ that can appear in the
process of applying the rules  $\^{\cal R}$ to $w$. To prove that 
$\cal T$ spans $\goth F$ it is enough to show the following two things:
\begin{itemize}
\item[(i)]
$w\not\in\Omega(w)$;
\item[(ii)]
$|\Omega(w)|<\infty$.
\end{itemize}
Indeed, if both (i) and (ii) are true, then the rules $\cal R$ can be
applied to $w$ at most finitely many 
times, after that $w$ is reduced to a linear combination of $\cal T$.

Recall that the set of generators $\cal B$ is linearly ordered. Extend
this order alphabetically to $\cal B^*$ by comparing the words from
left to right.
For a word $w=a_1(n_1)\cdots a_k(n_k)\in \cal X^*$ denote 
by $w_i=a_i(n_i)\cdots a_k(n_k)$ the $i$-th tail of $w$. Let also
$$
d(w) = -\sum_{i=1}^k n_i + \sum_{1\le i<j \le k} N(a_i,a_j).
$$
Take a word $u=b_1(p_1)\cdots b_k(p_k)\in \Omega(w)$. It is easy to
see that $(b_1, \ldots, b_k)$ is a permutation of $(a_1, \ldots, a_k)$.
Consider two sequences of integers: \break
$\big(d(w_1),\ldots, d(w_k)\big)$ and 
$\big(d(u_1),\ldots, d(u_k)\big)$. 
We note that  if $d(w)\ll 0$, then $w\in M(\cal Q)$, therefore
both claims (i) and (ii) follow from the following statement:

\begin{itemize}
\item[(iii)]{\sl
For every  $1\le i\le k$ we have 
$d(u_i)\le d(w_i)$; moreover, if\break $\big(d(w_1),\ldots, d(w_k)\big) =
\big(d(u_1),\ldots, d(u_k)\big)$, then 
$b_1\cdots b_k < a_1\cdots a_k\in \cal B^*$. }
\end{itemize}

It is enough to check (iii)
for a word $u$ that appears in the right-hand side of
a rule $w\ra f\in \cal R$. Then $u$ is either
$$
a_1(n_1)\cdots a_{j-1}(n_{j-1})
a_j(n_j-s)a_{j+1}(n_{j+1}+s)a_{j+2}(n_{j+2})\cdots a_k(n_k)
$$
for $s\ge 1$ or
$$
a_1(n_1)\cdots a_{j-1}(n_{j-1})
a_{j+1}(n_{j+1}+s)a_j(n_j-s)a_{j+2}(n_{j+2})\cdots a_k(n_k)
$$
for $s\le N=N(a_j,a_{j+1})$. In both cases 
$d(u_i) = d(w_i)$ when $i\neq j+1$.  
In the first case we always have $d(u_{j+1})=d(w_{j+1})-s<d(w_{j+1})$.
In the second case
\begin{align*}
d(u_{j+1}) &= d(w_{j+1}) -\sum_{i=j+2}^k N(a_{j+1},a_i)
+\sum_{i=j+2}^k N(a_j,a_i)+n_{j+1}-n_j+s \\
&\le d(w_{j+1}) +m_j - n_j + n_{j+1}
\end{align*}
where $m_j$ is given by \fl{m}. Now since $w\in M(\cal R)$, 
we have $n_j-n_{j+1} \le m_j$  and
the equality can hold only if $a_j > a_{j-1}$. This proves (iii).

\begin{Rem}
In fact we have shown that $\^{\cal R}$ is a rewriting system,
which is a generalization of the rewriting system constructed in
\cite{freecv} for the case when the locality bound is constant and
non-negative. Only minor modifications to the argument in
\cite{freecv} are needed to show that $\^{\cal R}$ is confluent, 
i.e. the final result of applications of the rules does not depend on
the order in which the rules are applied. Then the Diamond lemma
\cite{bergman,freecv} would imply that $\rho\cal T$ is a basis of
$\goth F$. However, in \sec{linind} we prove that $\rho\cal T$ is
linearly independent in $\goth F$ by a different method. 
\end{Rem}

\section{Proof of Theorem \lowercase{\ref{thm:basis}}, part II and proof of
Theorem \lowercase{\ref{thm:embedding}}}
\label{sec:linind}
Let $\goth F = \goth F_N(\cal B)$ be a free vertex algebra. 
Recall that there is a homomorphism $\f:\goth F \to V_\Lambda$ for 
$\Lambda = \Z[\cal B]$, see \sec{embedding}. Denote also $\Lambda_+ =
\Z_+[\cal B]$.
In \sec{span} we have proved that the image of the set 
$\cal T\in \cal X^*$
of words $a_1(n_1)\cdots a_k(n_k)\in \cal X^*$ satisfying the condition
\fl{nojumps} spans $\goth F$ over $\Bbbk$ under the projection
$\rho:\ol{\Bbbk\<\cal X\>}\to \goth F$. As before we will identify
$\cal T$ with $\rho\cal T$.
In this section we prove
that $\f\cal T\subset V_\Lambda$ is linearly independent over $\Bbbk$.
Combined with \sec{span} this will prove both \thm{basis} and \thm{embedding}. 

\paragraph{Step 1}
First of all we note that without a loss of generality we can assume that the form
$(\,\cdot|\,\cdot\,)$  is non-degenerate on $\Lambda$. 
For otherwise we can embed the
set $\cal B$ into a bigger set $\bar{\cal B}$ preserving the locality
bound $N$, such that the matrix $\{N(a,b)\}_{a,b\in\bar{\cal B}}$ is
non-degenerate. Then, if we have proved the statement for the
non-degenerate case, we know that the set $\cal T$ in linearly
independent in $V_{\Z[\bar{\cal B}\,]}$, hence
it is linearly independent in $V_\Lambda$.

\paragraph{Step 2}
Recall from \sec{confact} that for each function $f:\cal B\to \Bbbk$
there is a conformal derivation $\alpha_f\in cder\,\goth F$ such that 
$N(\alpha_f, \cal B) = 1$
and $\alpha_f(0)b = f(b)\,b$ for every  $b\in \cal B$ (see \lem{H+}\ref{H:f}).
Let $R=\Bbbk[\alpha_f(n)\,|\,f\in \cal B^*, \ n\in\Z_+]$ 
be the commutative polynomial algebra 
generated by all $\alpha_f(n)$. \lem{H+}\ref{H:h} states
that $R$ acts on both 
$\goth F$ and $V_\Lambda$ and these actions agree with $\f$.


Fix an element $a\in\cal B$ and let 
$f:\cal B\to \Bbbk$ be defined by $f(b)=\delta_{a,b}$. Denote the
corresponding conformal derivation by $a\check = \alpha_f \in cder\,
\goth F$. The algebra $R$ is graded by 
$\Lambda_+ \oplus \Z$ by setting $\op{wt}a\check(n) = a$, \ 
$\deg a\check(n) = -n$.

Our next observation reduces the problem to the following claim.

\smallskip\noindent
\bf Claim A.\quad \sl
For every homogeneous non-trivial linear combination  $x$ 
of the elements of $\cal T$  such that 
$\op{wt}x=\lambda$ and 
$\deg x > \deg w_{\min}(\lambda)=\frac12 (\lambda|\lambda)$,
there is $r\in R$, such that  
$\deg r = \deg w_{\min}(\lambda)- \deg x$
and $rx= c\, w_{\min}(\lambda)$ for some constant $0\neq c\in \Bbbk$. 
\rm\smallskip

\noindent
Here $w_{\min}(\lambda)$ is given by \fl{wmin}.

Indeed, since $\f$ is homogeneous, it is enough to show that
$\f(x)\neq 0$ for a homogeneous linear combination $x$ of $\cal T$. If
$x = w_{\min}(\lambda)$, then $\f(x)=\pm v_\lambda\neq0$ by
\prop{minimal}. Otherwise, using \lem{H+}\ref{H:h} we get
$r\f(x)=k\f\big(w_{\min}(\lambda)\big)\neq 0$, hence $\f(x)\neq0$.




\paragraph{Step 3}
Let $w=a_1(m_1)\cdots a_k(m_k)\in \cal X^*$ be a word in $\cal X$ and 
consider a monomial 
$p=a_1\check(n_1)\cdots a_k\check(n_k) \in R$. 
Let $G<S_k$ be a subgroup of the group of permutations of $k$ elements
consisting of permutations that fix the $k$-tuple 
$\left(a_1, \ldots, a_k\right)$, and let $H<G$ be the subgroup of $G$ 
consisting of permutations that in addition fix $(n_1, \ldots, n_k)$.  
We show that Claim A follows from

\smallskip\noindent
\bf Claim B.\quad \sl
There exists an element $r_p \in R$ such that 
\begin{equation}\label{fl:rpw}
r_pw =\frac1{|H|}
\sum_{\sigma\in G} a_1(m_1+n_{\sigma(1)})\cdots   
a_k(m_k+n_{\sigma(k)}).
\end{equation}
\rm\smallskip
Note that every word that appears in the right-hand side of \fl{rpw} has
coefficient 1. 

Indeed, suppose Claim B is true. Define an order on $\Z^k$ by letting
$\left(l_1,\ldots,l_k\right)>\left(l_1',\ldots,l_k'\right)$ 
if there is $1\le j\le k$ such
that $l_i = l_i'$ for $j+1\le i\le k$, but $l_j>l_j'$.  

Suppose that  $w\in \cal T$, and let
$(n_1,\ldots,n_k)=\eta(w)\in\Z^k$ be given by \fl{eta}.
Recall that $n_1\ge n_2 \ge \ldots\ge n_k$
and if $n_i=n_{i+1}$ then $a_i\le a_{i+1}$.
By \prop{minimal} we have $\rho\big(a_1(m_1+n_1)\cdots a_k(m_k+n_k)\big) =
w_{\min}(\op{wt}w)\neq 0$. Since $\cal T$ linearly spans $\goth F$, every
word of weight $\op{wt}w$ and degree less than $\deg w$ is 0 (in other
words, \prop{dong} holds). Therefore, if 
$\left(n_1',\ldots,n_k'\right)>\left(n_1,\ldots,n_k\right)$ then 
$\rho\big(a_1(m_1+n_1')\cdots a_k(m_k+n_k')\big)=0$. 

If $\sigma \in G\ssm H$, 
then $\left(n_1, \ldots, n_k\right) < 
(\,n_{\sigma(1)}, \ldots, n_{\sigma(k)}\,)$. 
Therefore, \fl{rpw} implies that
$$
r_p\, \rho(w) = \rho\big( 
a_1(m_1+n_1) \cdots a_k(m_k+n_k)\big)\neq 0.
$$
This proves Claim A for $x=w\in\cal T$.

Take now another word $u\in \cal T$ such that $u\neq w$, \ 
$\op{wt} u = \op{wt}w$, \ $\deg u = \deg w$. Let $\mu\in S_k$ be
a permutation such that 
$u = a_{\mu(1)}(m_1')\cdots a_{\mu(k)}(m_k')$. Set 
$(n_1',\ldots,n_k')=\eta(u)\in\Z^k$.
Suppose that $\left(n_1',\ldots,n_k'\right)\le 
\left(n_1,\ldots,n_k\right)$. Then we claim that $r_{p(w)}\, \rho(u)=0$. 
Indeed, we have 
$$
r_{p(w)}u = \sum_{\sigma\in \mu G\mu\inv} a_{\mu(1)}(m_1'+n_{\sigma(1)})\cdots   
a_{\mu(k)}(m_k'+n_{\sigma(k)}).
$$ 
This belongs to $\op{Ker}\rho$ if we show that 
$\left(n_1',\ldots,n_k'\right)<
\left(n_{\sigma(1)},\ldots,n_{\sigma(k)}\right)$
for all $\sigma \in \mu G\mu\inv$. 
If there were an equality in 
$\left(n_1',\ldots,n_k'\right)\le\left(n_1,\ldots,n_k\right)
\le \left(n_{\sigma(1)},\ldots,n_{\sigma(k)}\right)$,
then $\eta(u)$ would define the same colored partition of 
$\deg w-\deg_{\min}(\op{wt}w)$ as does $\eta(w)$, and this contradicts
the assumption $u\neq w\in\cal T$, see \sec{free}. 

Now we can prove Claim A. Let $x\in \Bbbk\<\cal X\>$ be a homogeneous
linear combination of words in $\cal T$. To every word 
$u = a_1(m_1)\cdots a_k(m_k) \in \cal T$ that appears in this
combination we correspond the sequence 
$\left(n_1,\ldots,n_k\right)=\eta(u)\in\Z^k$, given by \fl{eta}. 
All these sequences are pairwise different, since the words in $x$ are
all of the same weight and degree.
Let $w$ be the word that yields the maximal 
sequence. Take $p=p(w)\in R$ as above, and let $r_p$ be as in 
Claim B. Then the above argument shows that 
$0\neq r_p w \in \Bbbk w_{\min}(\op{wt} w)$ and $r_p u =0$ for every
other word $u\in\cal T$ that appears in $x$.  

\paragraph{Step 4}
It remains to construct for a word 
$w=a_1(m_1)\cdots a_k(m_k)\in \cal X^*$ and for a monomial
$p=a_1\check(n_1)\cdots a_k\check(n_k) \in R$ 
an element $r_p\in R$ such that \fl{rpw} holds.
Denote the right-hand side of \fl{rpw} by $S_p(w)$.

Let us introduce a partial order on the set of monomials in $R$.
Suppose  $a_i = a_j = a$ for some $1\le i\neq j\le k$, 
and consider the result of
substituting $a\check(n_i)a\check(n_j)$ by $a\check(n_i+n_j)a(0)$ in
$p$. We write $q\prec p$ if $q$ can be obtained by a number of such
substitutions. 

Let us now calculate the action of $p$ on $w$. Recall that $R$
acts on $\Bbbk\<\cal X\>$ by
derivations, so that
$a_i(n)\check\big(a_j(m)\big)=\delta_{i,j}\, a_j(m+n)$. 
For $1\le j\le k$ denote 
$\Theta_j = \set{i}{a_i=a_j}\subset \left\{1,\ldots,k\right\}$.
Then 
\begin{equation*}
pw = \sum_{i_1\in \Theta_1, \ldots, i_k\in \Theta_k}
a_1\Big(m_1+\sum_{j=1}^k \delta_{i_j,1}\, n_j\Big) 
\ \cdots \ a_k\Big(m_k+\sum_{j=1}^k \delta_{i_j,k}\, n_j\Big).
\end{equation*}
The crutial observation is that 
$$
pw=\sum_{q\preccurlyeq p} c(p,q)\, S_q(w)
$$
for some positive integers $c(p,q)$. Using this, we construct 
$r_p$ by induction on the number $l =l(p)= \#\set{1\le j\le k}{n_j\neq 0}$. 
If $l=0$ or 1, then $q\preccurlyeq p$ implies $q=p$ and we take
$r_p = c(p,p)\inv\, p$. Otherwise, set 
$$
r_p = c(p,p)\inv\Big(p - \sum_{q\precneqq p} c(p,q)\, r_q\Big).
$$

\hypertarget{biblio}{}
\bibliography{}

\end{article}

\end{document}